\documentclass[12pt,reqno,a4wide]{amsart}

\usepackage{tikz} 
\usepackage{calrsfs}
\usepackage[charter]{mathdesign}

\allowdisplaybreaks

\newtheorem{theorem}{Theorem}

\newtheorem{lemma}{Lemma}

\title{ How to sum powers of balancing numbers efficiently}

\author{Helmut Prodinger}
\address{Helmut Prodinger\\
	Department of Mathematical Sciences\\
	Stellenbosch University\\
	7602 Stellenbosch\\
	South Africa}
\email{hproding@sun.ac.za}

\keywords{Balancing numbers, Binet formula, generating functions}
\subjclass[2010]{11B49; 05A15}

\begin{document}
\begin{abstract}
	Balancing numbers possess, as Fibonacci numbers, a Binet formula. Using this, partial sums of arbitrary powers of balancing numbers can be summed explicitly. For this, as a first step, a power $B_n^l$ is expressed as a linear combination of $B_{mn}$.
	The summation of such expressions is then easy using generating functions.
		\end{abstract}

\maketitle
	
	\section{Introduction}
The balancing numbers $B_n$\footnote{There is no danger in this paper to confuse the notation with Bernoulli numbers or other popular sequences often denoted by $B_n$.} appear in a very recent paper \cite{Takao}. They are given by either the recurrence
\begin{equation*}
B_n=6B_{n-1}-B_{n-2},\quad B_0=0,\ B_1=1, 
\end{equation*}
or the generating function
\begin{equation*}
\sum_{n\ge0}B_nz^n=\frac{z}{1-6z+z^2};
\end{equation*}
they possess a Binet-type explicit formula
\begin{equation*}
B_n=\frac{\alpha^n-\beta^n}{4\sqrt2},
\end{equation*}
with
\begin{equation*}
\alpha=3+2\sqrt2,\quad \beta=3-2\sqrt2;
\end{equation*}
note that $\alpha\beta=1$ and $\alpha-\beta=4\sqrt2$. The sequence appears as A001109 in \cite{OEIS}; here are the first values:
\begin{equation*}
0, 1, 6, 35, 204, 1189, 6930, 40391, 235416, 1372105, 7997214,\dots
\end{equation*}

Apart from different $\alpha$, $\beta$, thanks to the Binet formula, the balancing numbers resemble Fibonacci numbers.

Recently, I have seen several papers (I don't want to provide explicit citations) where 
\begin{equation*}
\sum_{0\le k\le n}X_{km}^l
\end{equation*}
are computed, for, say, $X_n$ being Fibonacci numbers or other ones of a similar structure (one other example are Pell numbers).
Often, elementary methods are used that only work for a fixed small exponent. In \cite{Prodinger-Fro}, it was demonstrated how to deal with Fibonacci (and Lucas) numbers, and general exponents. Since it is felt that this approach should be more widely known, 
we work out the details of how to sum
\begin{equation*}
	\sum_{0\le k\le n}B_{km}^l
\end{equation*}
for balancing numbers; both $m, l$ are arbitrary positive integers.

\section{Linearization of powers of balancing numbers}

Our goal is to replace powers of balancing numbers by linear combinations of balancing numbers at different indices. The case of
odd exponents is simpler.

\begin{theorem}
For any $n\ge0$ and $l\ge0$
\begin{equation*} 
B_n^{2l+1}=\frac{1}{2^{5l}}\sum_{0\le s \le l}(-1)^s\binom{2l+1}{s}B_{(2(l-s)+1)n}.
\end{equation*}
\end{theorem}

\begin{theorem}
	For any $n\ge0$ and $l\ge1$
\begin{align*}
B_n^{2l}&=\frac{2}{2^{5l}}\sum_{0\le s <l}\frac{(-1)^s\binom{2l}{s}}{B_{2(l-s)}}
\Big(B_{2(l-s)n}+B_{2(l-s)(n+1)}\Big)\\&
+\frac{1}{2^{5l}}\sum_{0\le s <l}\frac{(-1)^{s-1}\binom{2l}{s}B_{2(l-s)}}{B_{l-s}^2}B_{2(l-s)n}
+(-1)^l\binom{2l}{l}.
\end{align*}
\end{theorem}

Essentially, these identities are polynomial identities, if we write $\alpha=x$, $\beta=\frac1x$. By polynomial we allow here also negative exponents, but there are only finitely many terms. We prove the easier case of odd exponents; the other one is similar but longer. The quadratic equation that $\alpha$ resp.\ $\beta$ satisfy is only used in the more involved case of even exponents.

\begin{align*}
B_n^{2l+1}&=\frac{(x^n-x^{-n})^{2l+1}}{2^{5l}(x-\frac1x)}\\
&=\frac1{2^{5l}(x-\frac1x)}{\sum_{s=0}^{2l+1}(-1)^s\binom{2l+1}{s}x^{sn-(2l+1-s)n}}\\
&=\frac1{2^{5l}(x-\frac1x)}{\sum_{s=0}^{l}(-1)^{s}\binom{2l+1}{s}\Big(x^{(2l+1-2s)n}-x^{-(2l+1-2s)n} \Big)}\\
&=\frac1{2^{5l}}{\sum_{s=0}^{l}(-1)^{s}\binom{2l+1}{s}B_{(2l+1-2s)n}}.
\end{align*}

\section{Summing equally spaced balancing numbers}

Now we deal with the sum
\begin{equation*}
\sum_{0\le k \le n}B_{km}=[z^n]\frac1{1-z}\sum_{k\ge0}B_{km}z^k
\end{equation*}
for any positive integer parameter $m$. (The notion $[z^n]f(z)$ refers to the coefficient of $z^n$ in the power series expansion of $f(z)$, as is now customary.) Once we understand this, we can apply such a formula to each term of the righthand sides
of the previous section, and our goal has been achieved.

\begin{lemma}
\begin{equation*}
	\sum_{k\ge0}B_{km}z^k=\frac{B_mz}{1-(6B_m-2B_{m-1})z+z^2}.
\end{equation*}
\end{lemma}
\textbf{Proof.}
The equivalent statement is
\begin{equation*}
	(1-(6B_m-2B_{m-1})z+z^2)\sum_{k\ge0}B_{km}z^k=B_mz,
\end{equation*}
which, in terms of coefficients, means for $m=1$
\begin{equation*}
	[z^1](1-6z+z^2)\sum_{k\ge0}B_{k}z^k=1,
\end{equation*}
which is correct; now we assume that $m\ge2$ and need to prove
\begin{equation*}
B_{km}-6B_mB_{(k-1)m}+2B_{m-1}B_{(k-1)m}+B_{(k-2)m}=0. 
\end{equation*}
In terms of the Binet form, this means
\begin{multline*}
(\alpha^{km}-\beta^{km})(\alpha-\beta)-6(\alpha^{m}-\beta^{m})(\alpha^{{(k-1)}m}-\beta^{{(k-1)}m})\\
+2(\alpha^{m-1}-\beta^{m-1})(\alpha^{{(k-1)}m}-\beta^{{(k-1)}m})+(\alpha^{{(k-2})m}-\beta^{{(k-2)}m})(\alpha-\beta)=0.
\end{multline*}
It can be checked (best with the help of a computer) that this is indeed true. \qed

\begin{theorem} For all $n\ge0$, $m\ge1$, we have
	\begin{equation*}
\sum_{0\le k\le n}B_{km}=\frac{B_{m(n+1)}-B_{mn}-B_m}{-2+6B_m-2B_{m-1}}.
	\end{equation*}
\end{theorem}

\textbf{Proof.} This follows from partial fraction decomposition
\begin{multline*}
\frac1{1-z}\sum_{k\ge0}B_{km}z^k	
=\frac{B_m}{-2+6B_m-2B_{m-1}}\frac{1-z}{1-(6B_m-2B_{m-1})z+z^2}-\frac{B_m}{-2+6B_m-2B_{m-1}}\frac1{1-z}
\end{multline*}
and extraction of the coefficient of $z^n$ on both sides.\qed

\section{Conclusion}
This general method (presented on 4 pages) instead of a special case (presented on 16 pages, as recently seen) is by no means restricted to balancing number and/or Fibonacci numbers. It is essential to have a Binet-type formula. Then sums, alternating sums, weighted sums like $\sum_{0\le k\le n}kB_k$, etc. can all be computed explicitly, using generating functions. To use a computer algebra system is essential, as one is often forced to first guess a formula before it is then proved.

	\bibliographystyle{plain}

\end{document}